\documentclass[a4paper,12pt,twoside]{article}

\usepackage[body={165mm,235mm}]{geometry}
\usepackage{a1}
\usepackage[english]{babel}
\usepackage[latin1]{inputenc} 
\usepackage{amsfonts,amssymb}
\usepackage{amsmath}
\usepackage{graphicx}	

\usepackage{makeidx}
\makeindex

\def\mapright#1#2#3{\smash{\mathop{\hbox to
#3{\rightarrowfill}}\limits^{#1}_{#2}}}

\def\mapleft#1#2#3{\smash{\mathop{\hbox to
#3{\leftarrowfill}}\limits^{#1}_{#2}}}

\def\mapright#1#2{\smash{\mathop{\hbox to 0.90cm{\rightarrowfill}}\limits^{#1}_{#2}}}
\def\mapleft#1#2{\smash{\mathop{\hbox to 0.90cm{\leftarrowfill}}\limits^{#1}_{#2}}}

\def\mapleftright#1#2{\smash{\mathop{\hbox to 0.80cm{\leftarrowfill \rightarrowfill}}\limits^{#1}_{#2}}}

\title{A tougher challenge to 3-manifold topologists \\ and  group algebraists
} 
\author{Sóstenes L. Lins}

\date{\today}

\begin{document}

\maketitle

\begin{abstract}
This paper poses some basic questions about instances (hard to find) of 
a special problem in 3-manifold topology.
``Important though the general concepts and propositions may be with the 
modern industrious passion
for axiomatizing and generalizing has presented us \ldots nevertheless I am convinced that
the special problems in all their complexity constitute the stock and the core of mathematics; 
and to master their difficulty requires on the whole the harder labor.'' Hermann Weyl 1885-1955, 
cited in the preface of the first edition (1939) of A. N. Whitehead's book
{\em The classical groups: their invariants and representations} \cite{whitehead1997}.

In this paper I focus on new uncertainties left unanswered in L. Lins thesis \cite{lins2007blink}
on the homemorphism problem of eleven concrete pairs of closed orientable 3-manifolds induced
by 3-connected monochromatic {\em blinks} (\cite{kauffman1994tlr}). The eleven HG8QI-classes are the 
only doubts left in the thesis, but the first two of them were solved few days ago and 
in this work I report on their solutions. We also include an appendix which can be used
to import all the links of this paper into SnapPy. The appendix was obtained by drawing the links
in SnapPy, work performed by C. Nascimento.

\end{abstract}

\section{Introduction}
In a joint recent paper posted recently in the arXiv (\cite{linslinschallenge2013}) 
my son Lauro Lins and myself ask some 6 years old questions for which we had no answers
about homeomorphisms between closed orientable 3-manifolds. The two pairs of 3-manifolds
were the only uncertainties that were left in L. Lins thesis  
(\cite{lins2007blink}) under my supervision in the 
domain of 3-manifolds being induced by arbitrary connected blinks
up to 9 edges (9-small 3-manifolds).
A subset of relevant 10-crossings blinks were generated but their topological classification remains untouched.
The paper was taken seriously by a few researchers, among them M. Culler, N. Dunfield, C. Hodgson 
and others that could solve them very quickly using GAP (\cite{gap2002gap}), Sage (\cite{sage2012}) 
and SnapPy (\cite{snappy}), 
tools that (except for GAP) were basically unknown to us. 
The solutions were obtained by distinct methods and are all consistent (inclusive with BLINK, 
the program of L. Lins (implementing my theory described in \cite{lins1995gca}), which support his thesis). Together
with my colleague Cristiana Nascimento, here at CIn/UFPE, I am learning fast to operate these wonderful tools.
The solutions people found shows that BLINK does a complete job in
topologically classifying the 9-small 3-manifolds. 
This is the subject of a joint paper with Lauro, currently under preparation.

The first solution that I got, and that still blows my mind,  
was by Craig Hodgson using length spectra techniques, based in his
joint paper with J. Weeks entitled {\em Symmetries, 
isometries and length spectra of closed hyperbolic three-manifolds} 
(\cite{hodgson1994symmetries}). By using SnapPy Craig showed that even though the 
quantum WRT-invariants as well as the volumes of the 
hyperbolic $Z$-homology spheres induced by the bfl's,
$U[1466]$ and $U[1563]$ are the same, the length of the 
smallest geodesics of them are distinct. For the other pair of bfl's $U[2125]$ and $U[2165]$ he shows that 
precisely the same facts apply. Here is a summary of Craig's findings extracted from the SnapPy session
that he kindly sent me. As Craig writes: {\em ``The output of the length spectrum command shows the {\em complex lengths}
of closed geodesics --- the real part is the actual length and the imaginary part is the rotation angle
as you go once around the geodesic.''}\\
\begin{center}
\center{Class $9_{126}$:}
\begin{verbatim}
First geodesic of U[1466]: 1.0152103824828331+0.39992347315914334i.
First geodesic of U[1563]: 0.9359206605025168+2.333526236965665i.
Volume of both manifolds: 7.36429600733.
\end{verbatim}
\center{Class $9_{199}$:}
\begin{verbatim}
First geodesic of U[2125]:  0.8939075859248593+0.761197185679321i.
First geodesic of U[2165]:  0.7978548001747316+2.9487425029345973i.
Volume of both manifolds: 7.12868652133.
\end{verbatim}
\end{center}

I posted 4 versions of \cite{linslinschallenge2013} correcting annoying mistakes 
in the presentations of the fundamental groups, putting a second pair of links,  
and in focusing the challenge in a broader context.  
I computed the presentations manually and I had a hard time making them correct.
Even though the presentations are redundant because the blink is enough to 
define the 3-manifold, as explained in  \cite{kauffman1994tlr},
my objective was to facilitate the work for those wanting to use GAP.
The time spanned between the first (April 22, 2013) and the last version (May 1, 2013) 
was a little more than one week. 
During these revisions I was completely unaware that the paper had called the attention
of many people. I did not know that the blog on lower dimensional topology
was very active exposing my incorrections and I apologize for my ignorance. 
I thank Cristiana for having calling the blog to my attention. Worse, 
some people did not see the follow up versions.
This was the case of Nathan Dunfield who worked with the 
wrong presentations. Not without reason he was angry at me,
but I think that this is no longer true, since he was willing to answer my sometimes naives and stupid 
questions and send me a solution for the first pair of manifolds of the present work, 
using SnapPy, Sage and GAP computations, by working with covers. I did not know these tools.
But, when properly motivated, I can learn fast and in general I do believe 
that I have something important and different to say in this brave new world of
3-manifolds: see the wonderful essay of E. Klarreich published by the Simons Foundation (march 2012), 
\cite{Klarreich1012}. I have been putting a great amount of time and effort during my scientific carreer, 
(most of the time as an isolated researcher) on (mainly closed) 3-manifods. I seek no longer to be isolated: my team
is the World, my compromisse is with Truth (independently of whom first found it).

Marc Culler was very helpful in answering questions of myself and Cristiana and helping her about issues in the 
downloading and 
installing SnapPy and Sage and GAP in her machine. With the presentation incorrections out of the 
way he produced an independent proof of the distinctveness of $(U[1466, U[1563])$ and of 
$(U[2125, U[2165])$. He also produced instantaneous isomorphic triangulations of the homeomorphic
3-manifolds in the classes $9_{126}$ and $9_{199}$. This fact makes me anxious to compare and timing the 
performances of BLINK (which also produces instantaneous solutions
for the same problems) and SnapPy regarding 
finding homeomorphisms of $k$-small 3-manifolds, given that the 
homeomorphisms exist.

\section{Objective of this work: help to make BLINK known}
In this paper I put some new challenges (also coming from
\cite{lins2007blink}), that seem harder than the ones considered in the 
previous paper. The reason I think so
is that going from 9 to 14,15,16 crossings in the links,
numerical problems start appearing concerning finding the Dirichlet domain and, in these cases, finding
isomorphic triangulations might be harder to SnapPy than to BLINK. At any rate I have hundreds of examples
where the performance of these programs in this issue could be compared, if anyone is seriously interested. Currently
BLINK is not documented and one of my objectives is to seek for help in doing it and extend its capability.
BLINK is hosted at Github under the userid {\em laurolins} and is open source code project. Unfortunately Lauro (currently
a researcher at AT\&T) does not have the necessary time to go on with the implementation. But he welcomes and 
is willing to help collaborators in getting started. 
As for myself, I am too old for the energy needed to construct good pieces of software. I intend to act
as one of some Scientific Supervisors for the deployment and for the discussions of the new algorithms  
to be included in BLINK, but only at the mathematical level. The technological and software engineering
screws and bolts needed, I leave to others.

An algorithm that I want to attach to BLINK is finding a uniformly distributed random closed orientable 3-manifold induced 
by a blink with an arbitrary number (even thousands) of edges. I want to gather evidence 
for the truth of some important conjectures that depend on this
capability. Another example of such new algorithms that I want to include in BLINK is made possible by the theory 
in Ricardo Machado's thesis
under my supervision, defended in March, 2013. We got an $O(n^2)$-algorithm for going from a special kind of gem, 
named {\em resoluble gem}, to 
a blink inducing the same manifold. This work is available, in still rather sketchy form (even the definition of 
resolubility is unecessary complicated), in the three joint papers
posted last year in the arXiv, {\cite {linsmachadoA2012, linsmachadoB2012, linsmachadoC2012}. The algorithm was
implemented in Mathematica, but it needs to be improved and re-implemented in Java or C++. 
We found a rather crude framed link presentation for
the hyperbolic dodecahedral space (Weber-Seifert manifold). As far as I know nobody has found such a
framed link. My interest in it was aroused by J. Weeks in a visit to the Geometry Center in April 1993, when he asked me whether
I had such framed link. The link inducing the Weber-Seifert 3-manifold is a 9-component link embedded 
into $\mathbb{R}^3$, with an integer attached to each component (its framing) 
and having a total of 68 (only) vertices with a projection having 142 crossings.
(It started with a PL-link with more than 600 vertices.)
In a fourth joint paper with R. Machado, currently under preparation,
we will show that every 3-manifold admits a resoluble gem inducing it.

\section{The 11 $HG8QI_t$-classes of blinks left unresolved in \cite{lins2007blink}}
I assume that the reader has with him a copy of the version 4 of previous challenge paper 
(\cite{linslinschallenge2013}) 
and has learned how to read the manifold either from the blink  or
from the blackboard framed link, \cite{kauffman1994tlr}. As for obtaining a presentation
of the fundamental group based on the Wirtinger relators  (\cite{stillwell1993classical}) 
and the Dehn fillings (\cite{rolfsen2003knots})
the two detailed examples given in \cite{linslinschallenge2013} should suffice, if the reader
has not available other pieces of softwares to get the presentation by automatic means. 
Actually, the best way to enter these manifolds
is to draw the blackboard framed link using SnapPy and informing the $w$ of each component as
its self-writhe in the projection. The framing of that component to be informed to SnapPy is $(w,1)$.
The complex numbers in polar form which appear at each $m^t_{p}$-class are the common 
quantum WRT-invariants. All except one of the eleven classes are formed by $Z$-homology spheres.
The exception is $16^t_{56}$ which has no torsion but Betti number 1. These facts are indicated 
by the small number in parenthesis (which  gives the homology of the manifold).
Actually the first two classes were recently resolved, only remaining the nine final ones.
Here is N. Dunfield's Sage session distinguishing the two manifolds 
induced by the two blinks in $14^t_{24}$:
\begin{verbatim}
sage: from snappy import *
sage: M1 = Manifold('1424_T71.tri')
sage: M2 = Manifold('1424_T79.tri')
sage: covers1 = M1.covers(5, method='gap')
sage: covers2 = M2.covers(5, method='gap')
sage: [C.homology() for C in covers1]
[Z/132 + Z/132, Z/63 + Z/63, Z/3 + Z/3 + Z/3 + Z/3]
sage: [C.homology() for C in covers2]
[Z/3 + Z/3 + Z/3 + Z/3, Z/213 + Z/213, Z/432 + Z/432]
\end{verbatim}
and Cristiana's Sage session also 
distinguishing the two manifolds induced by the two blinks in $14^t_{24}$ and, 
in conjuction with BLINK, topologically classifying the manifolds induced by  the four blinks in $15^t_{16}$:
\begin{verbatim}
M=1424_T71, N=1424_T79
sage: [C.homology() for C in coversM]
[Z/3 + Z/3 + Z/3 + Z/3, Z/63 + Z/63, Z/132 + Z/132]
sage: [C.homology() for C in coversN]
[Z/3 + Z/3 + Z/3 + Z/3, Z/213 + Z/213, Z/432 + Z/432]
------------------------------------------------------
A=1516_T118, B=1516_T119, C=1516_T181, D=1516_T205
sage: [X.homology() for X in coversA]
[Z/229773, Z/1110327, Z/3699687, Z/3018207]
sage: [X.homology() for X in coversC]
[Z/1110327, Z/229773, Z/3018207, Z/3699687]
sage: [X.homology() for X in coversB]
[Z/1052067, Z/3 + Z/1299909, Z/4117827, Z/126627]
sage: [X.homology() for X in coversD]
[Z/4117827, Z/1052067, Z/3 + Z/1299909, Z/126627]
\end{verbatim}

\eject
\subsection{The $HG8QI_t$ class $14^t_{24}$:}
\begin{figure}[!h]
\begin{center}
\includegraphics[width=10.0cm]{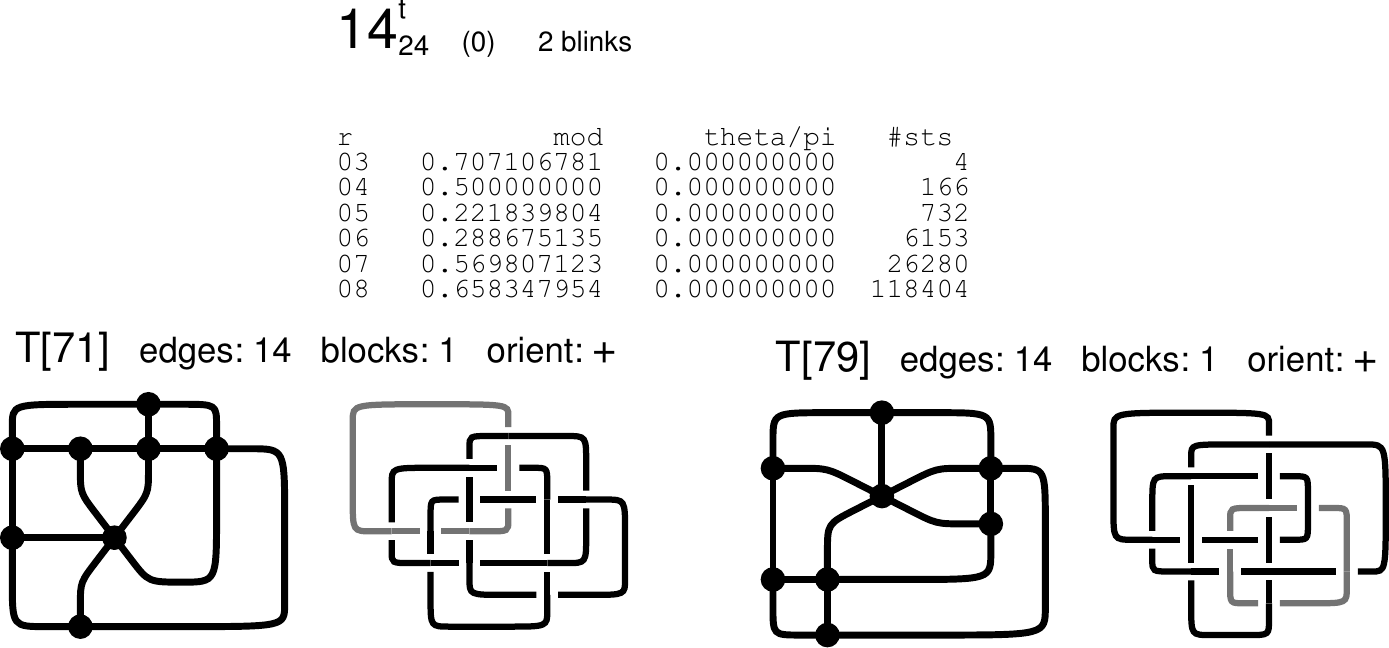}
\caption{\sf The above two manifolds are not homeomorphic. They are distinguished by the homology
of their 5-covers. This was immediately noted by N. Dunfield using Sage and GAP from triangulations
obtained by C. Nascimento using SnapPy, which could not find the Dirichlet domain due to numerical instability. 
}
\label{fig:firstdoubtA}
\end{center}
\end{figure}
\subsection{The $HG8QI_t$ class $15^t_{16}$:}
\begin{figure}[!h]
\begin{center}
\includegraphics[width=10.0cm]{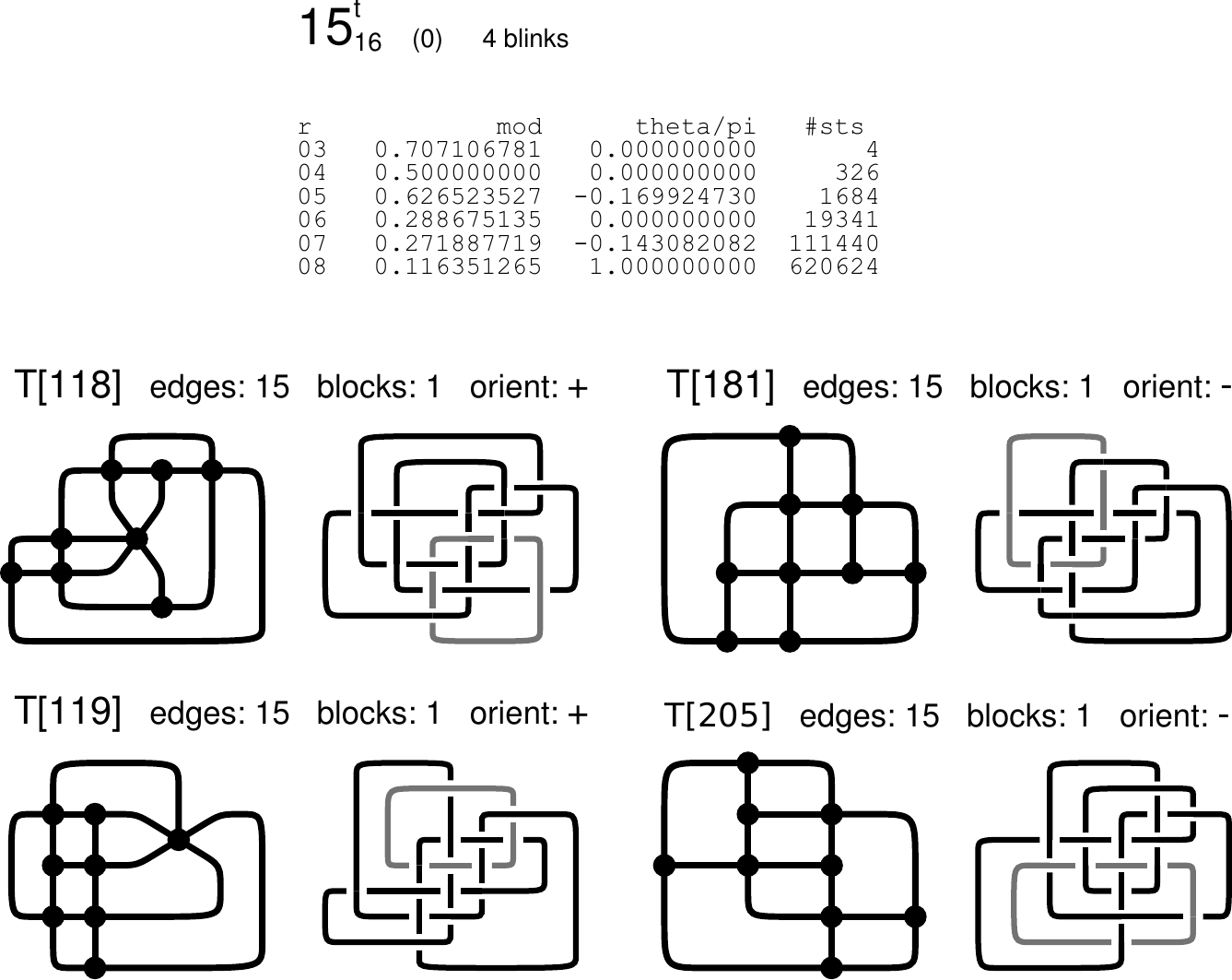}
\caption{\sf The above two manifolds are also non-homeomorphic. 
They are also distinguished by the homology
of their 5-covers. Relative to the class $15^t_{24}$ class $15^t_{16}$
the Sage/GAP software demands much more time.
This was obtained by C. Nascimento using SnapPy/Sage/GAP. The software
SnapPy could not find the Dirichlet domain due to numerical instability. 
}
\label{fig:T15_16}
\end{center}
\end{figure}

\eject
\subsection{The $HG8QI_t$ class $15^t_{19}$:}
\begin{figure}[!h]
\begin{center}
\includegraphics[width=9cm]{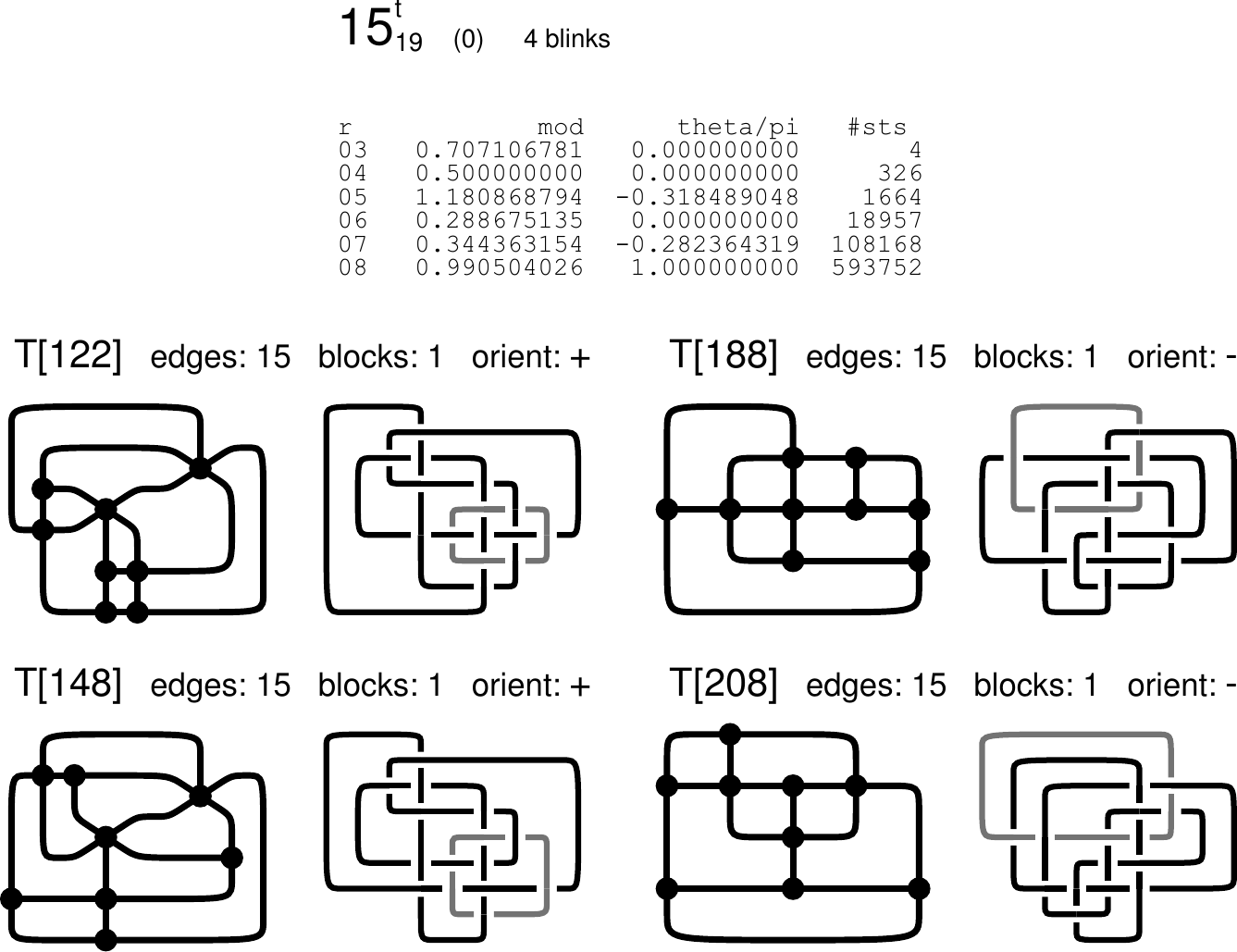}
\caption{\sf I do not know whether the above four manifolds are homeomorphic or not. 
BLINK says that there are at most two homeomorphisms classes among the four and I bet that 
this bound is attained.}
\label{fig:T15_19}
\end{center}
\end{figure}
\subsection{The $HG8QI_t$ class $15^t_{22}$:}
\begin{figure}[!h]
\begin{center}
\includegraphics[width=9cm]{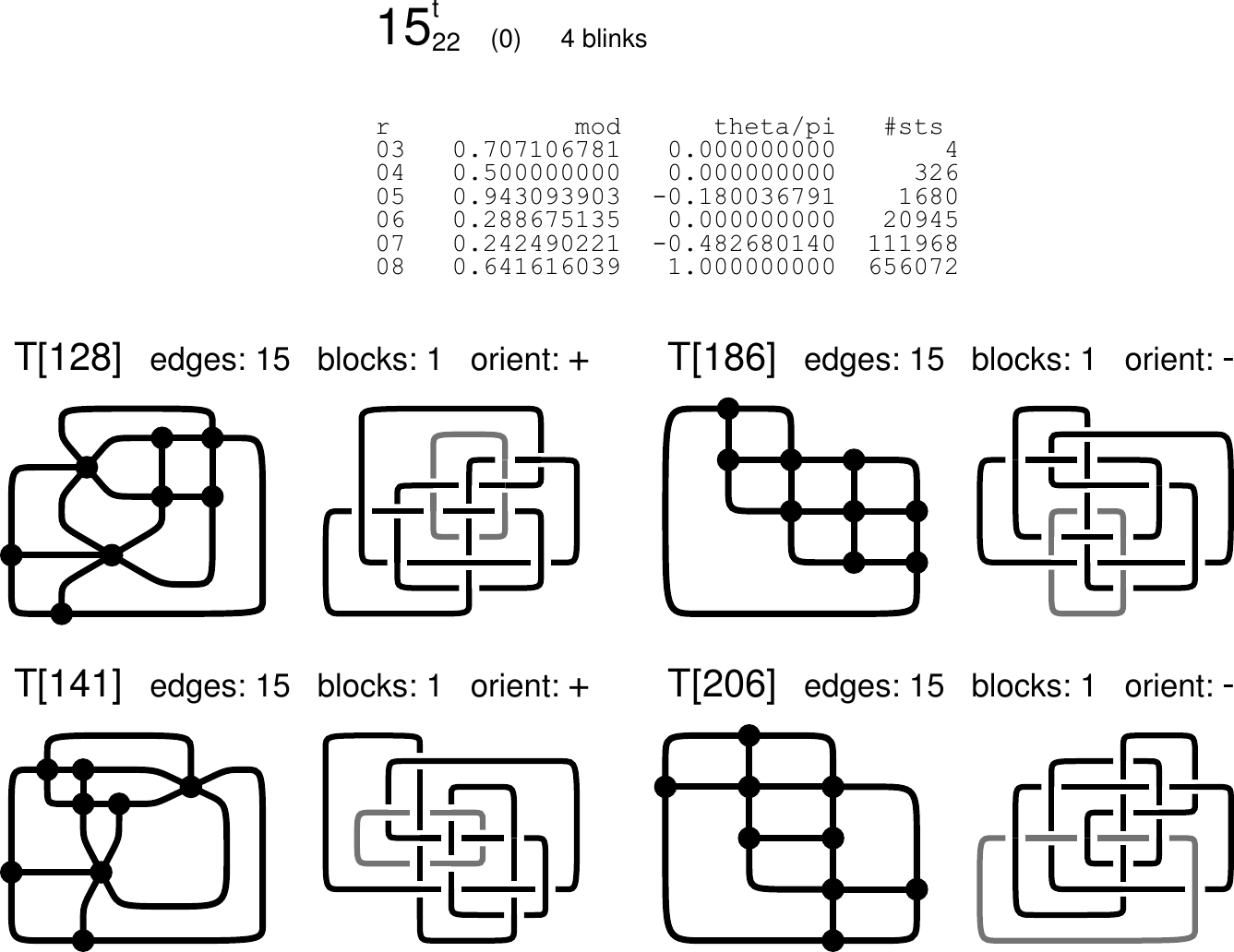}
\caption{\sf I do not know whether the above four manifolds are homeomorphic or not. 
BLINK says that there are at most two homeomorphisms classes among the four and I bet that 
this bound is attained.}
\label{fig:T15_22}
\end{center}
\end{figure}

\eject
\subsection{The $HG8QI_t$ class $16^t_{42}$:}
\begin{figure}[!h]
\begin{center}
\includegraphics[width=9cm]{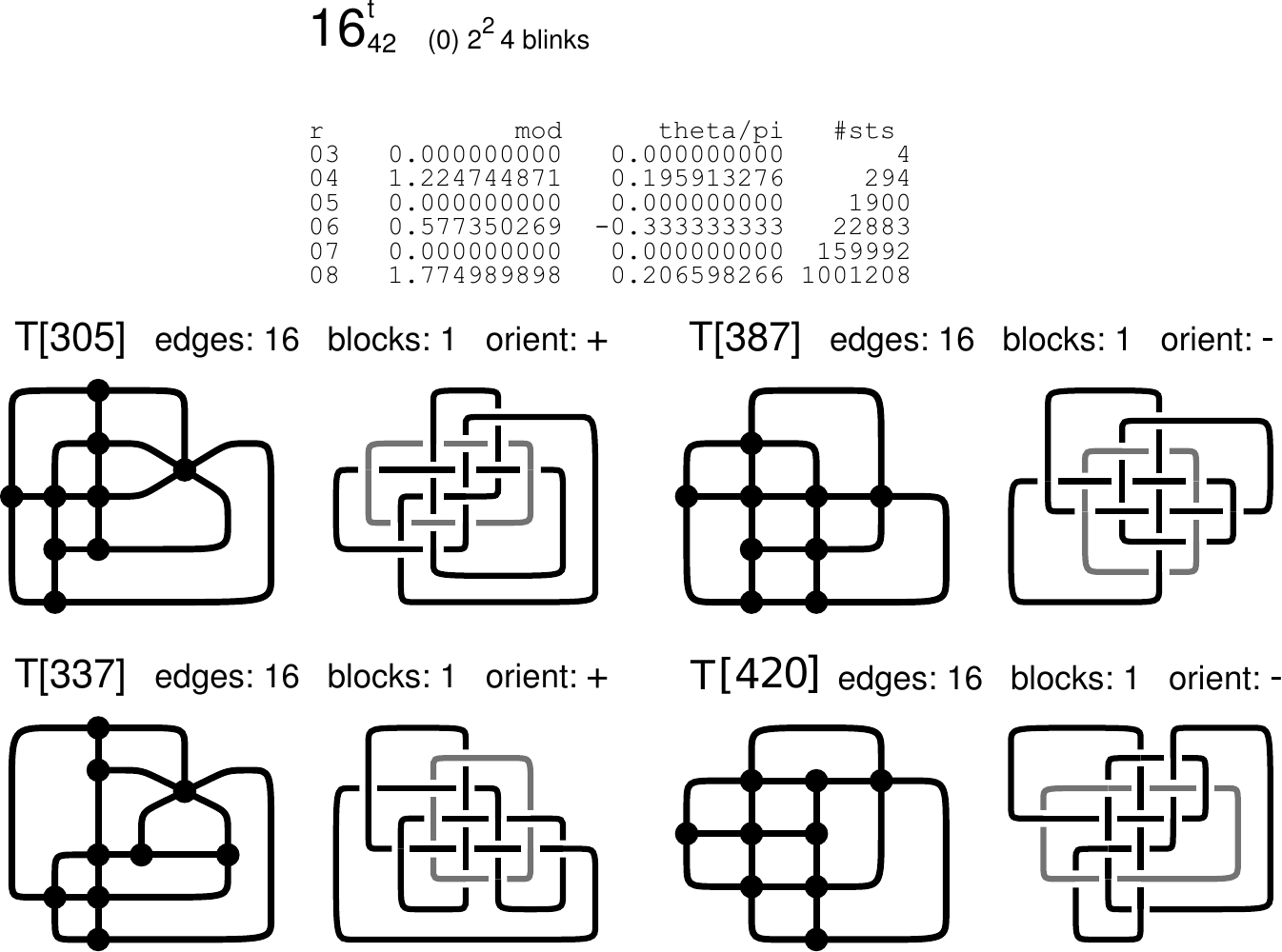}
\caption{\sf I do not know whether the above four manifolds are homeomorphic or not. 
BLINK says that there are at most two homeomorphisms classes among the four and I bet that 
this bound is attained.}
\label{fig:T16_42}
\end{center}
\end{figure}
\subsection{The $HG8QI_t$ class $16^t_{56}$:}
\begin{figure}[!h]
\begin{center}
\includegraphics[width=9cm]{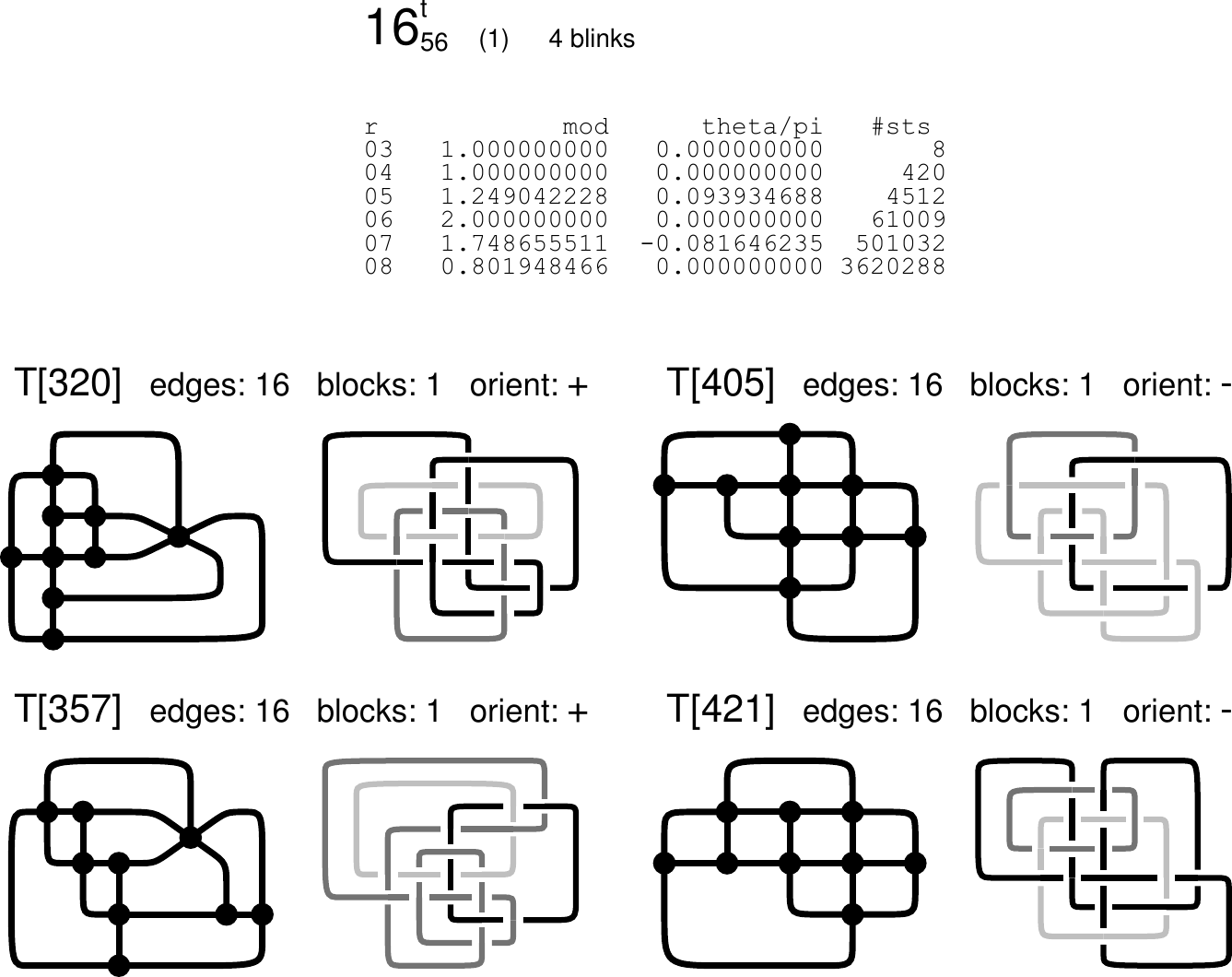}
\caption{\sf I do not know whether the above four manifolds are homeomorphic or not. 
BLINK says that there are at most two homeomorphisms classes among the four and I bet that 
this bound is attained.}
\label{fig:T16_56}
\end{center}
\end{figure}

\eject
\subsection{The $HG8QI_t$ class $16^t_{140}$:}
\begin{figure}[!h]
\begin{center}
\includegraphics[width=11cm]{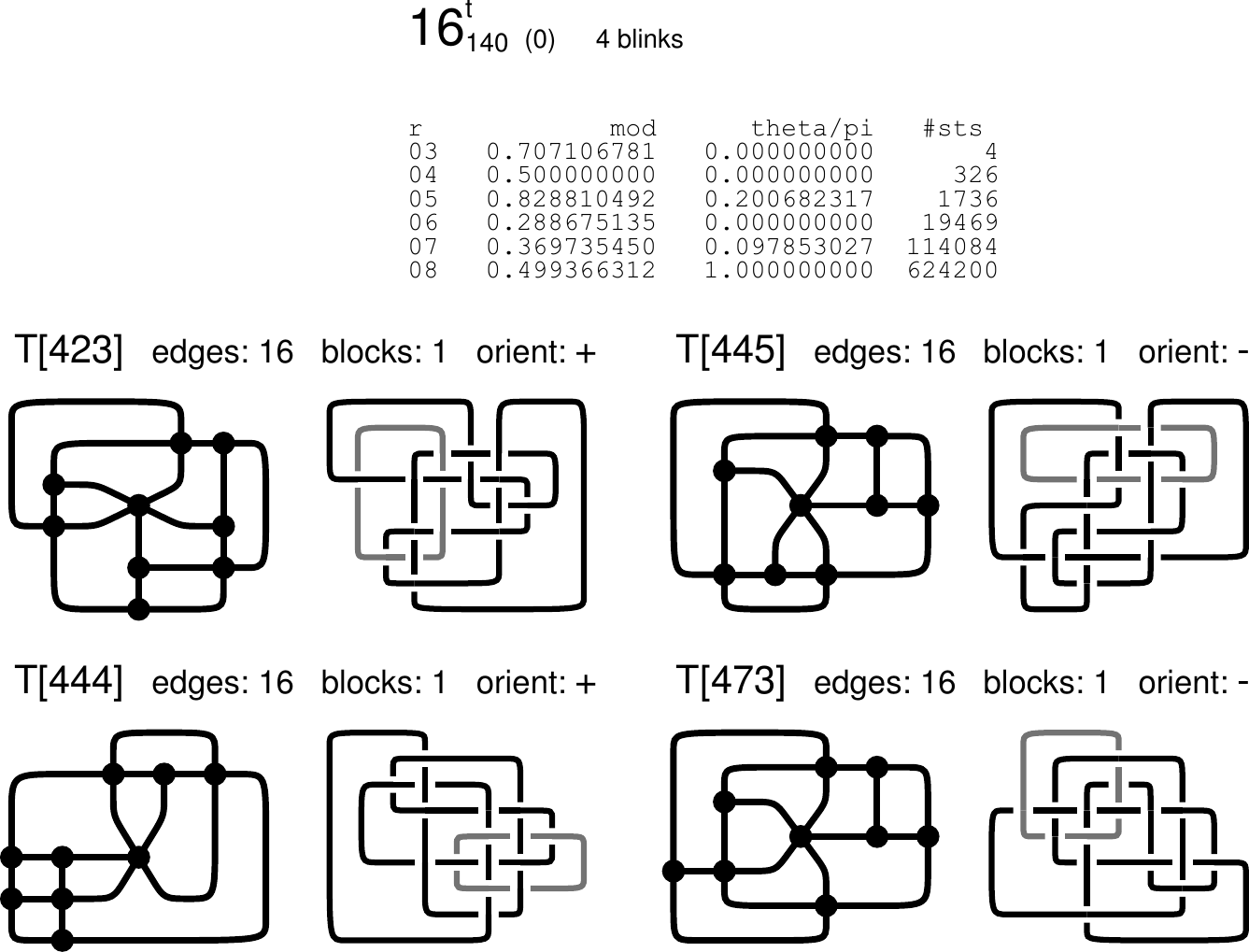}
\caption{\sf I do not know whether the above four manifolds are homeomorphic or not. 
BLINK says that there are at most two homeomorphisms classes among the four and I bet that 
this bound is attained.}
\label{fig:T16_140}
\end{center}
\end{figure}
\subsection{The $HG8QI_t$ class $16^t_{141}$:}
\begin{figure}[!h]
\begin{center}
\includegraphics[width=11cm]{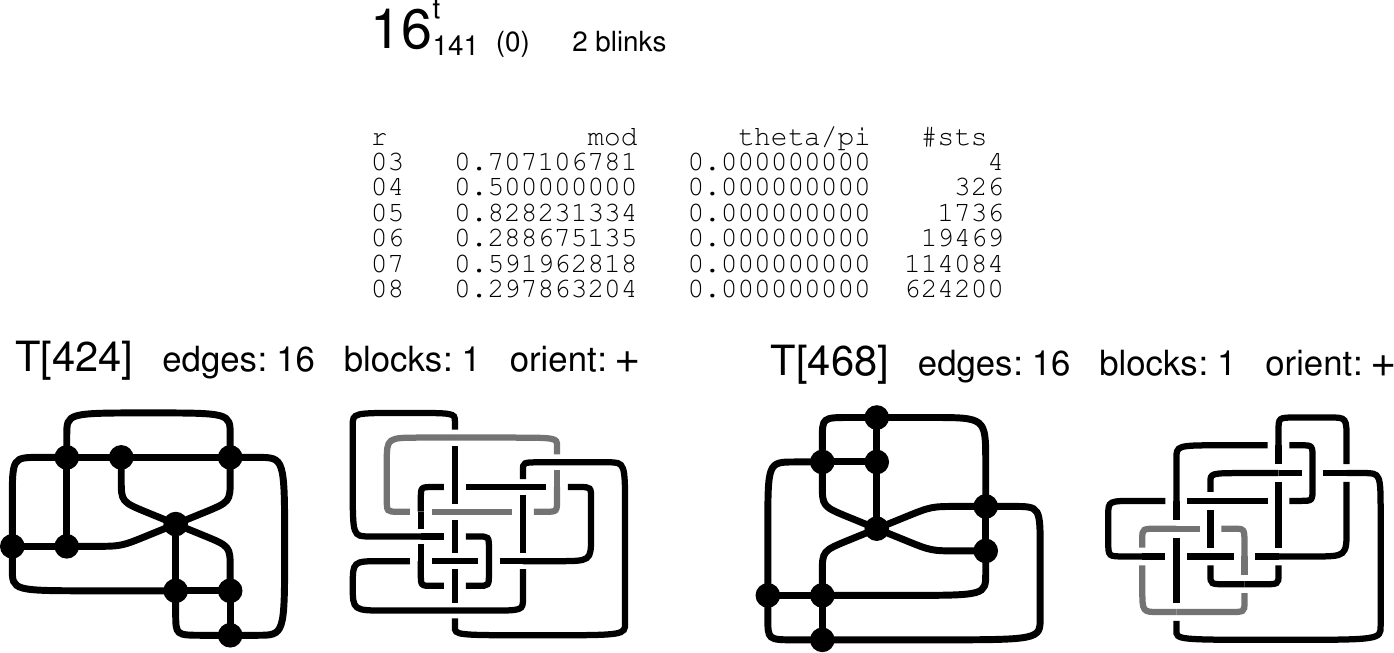}
\caption{\sf I do not know whether the above two manifolds are homeomorphic or not.}
\label{fig:T16_141}
\end{center}
\end{figure}

\eject
\subsection{The $HG8QI_t$ class $16^t_{142}$:}
\begin{figure}[!h]
\begin{center}
\includegraphics[width=12cm]{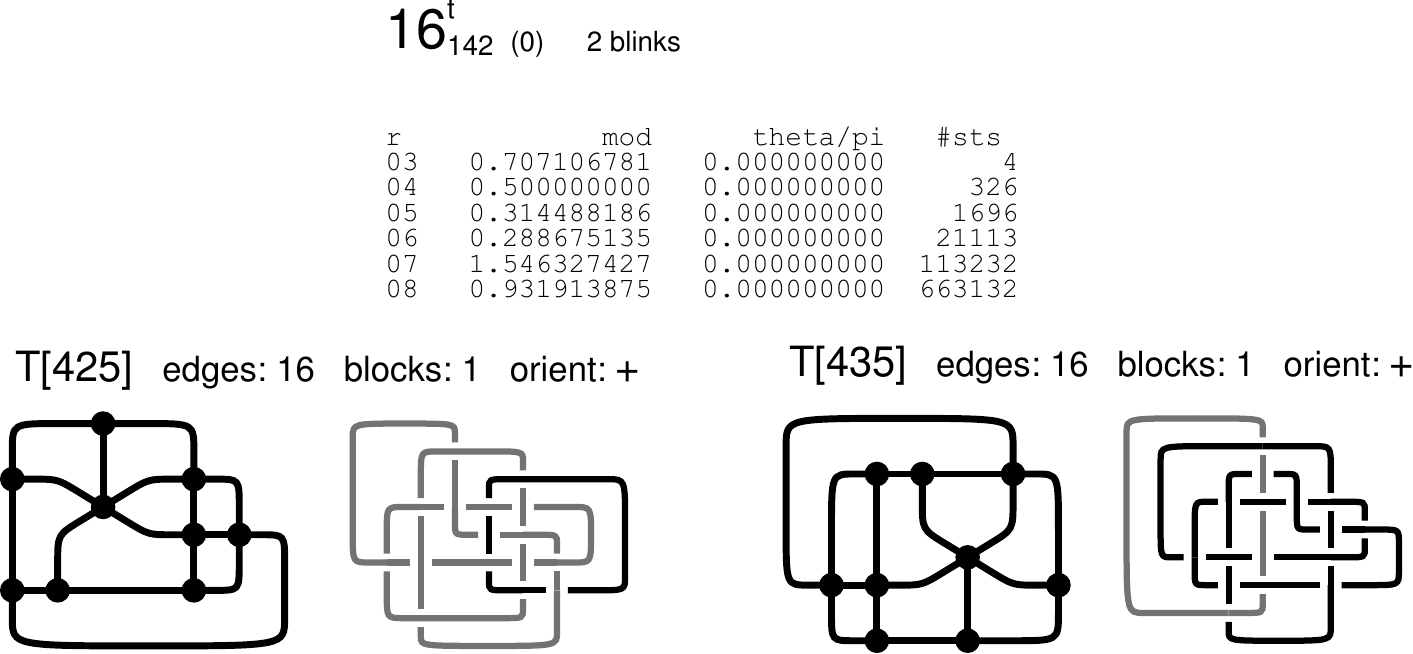}
\caption{\sf I do not know whether the above two manifolds are homeomorphic or not.}
\label{fig:T16_142}
\end{center}
\end{figure}
\subsection{The $HG8QI_t$ class $16^t_{149}$:}
\begin{figure}[!h]
\begin{center}
\includegraphics[width=12cm]{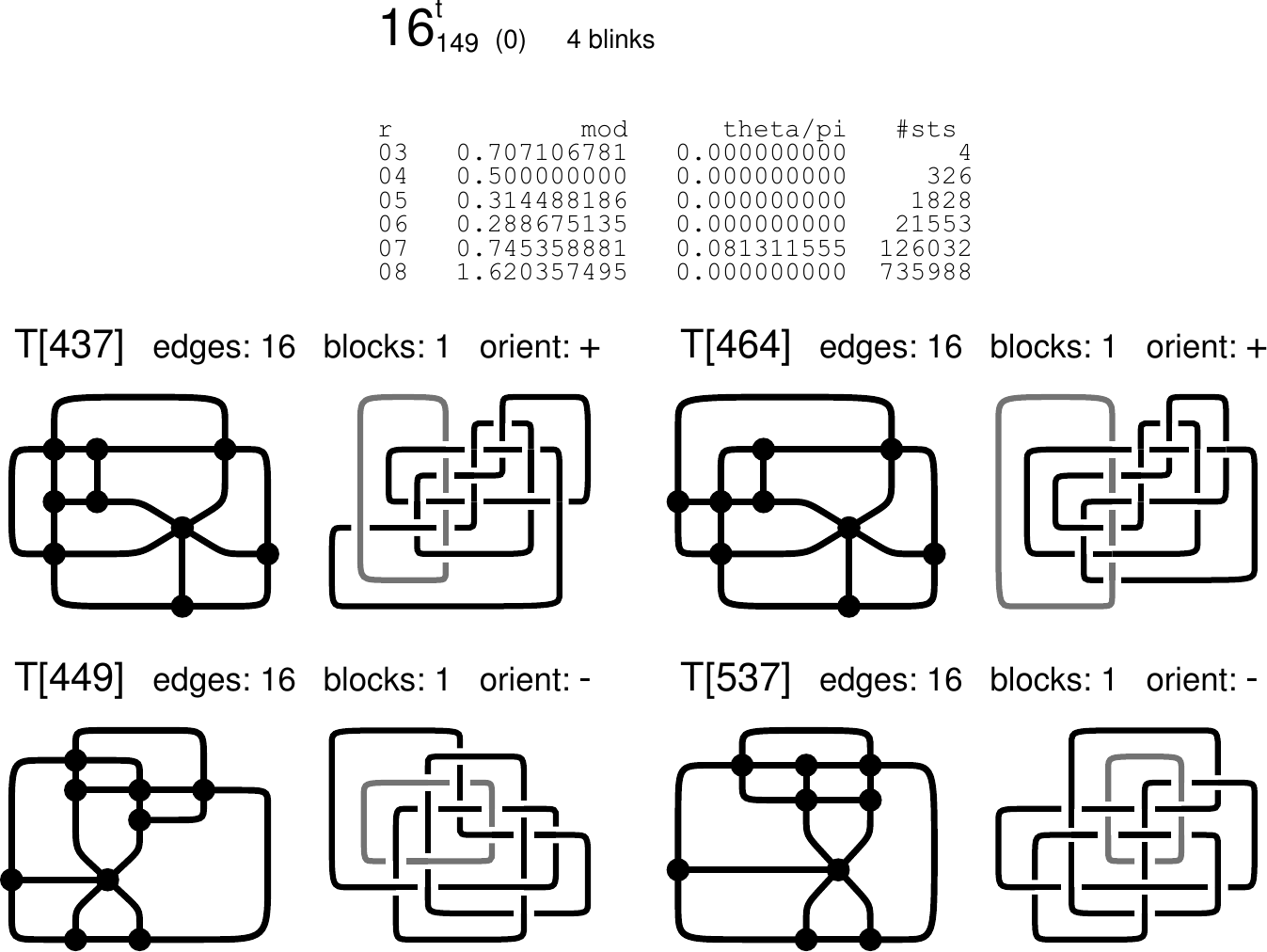}
\caption{\sf I do not know whether the above four manifolds are homeomorphic or not. 
BLINK says that there are at most two homeomorphisms classes among the four and I bet that 
this bound is attained.}
\label{fig:T16_149}
\end{center}
\end{figure}

\eject
\subsection{The $HG8QI_t$ class $16^t_{233}$:}
\begin{figure}[!h]
\begin{center}
\includegraphics[width=12cm]{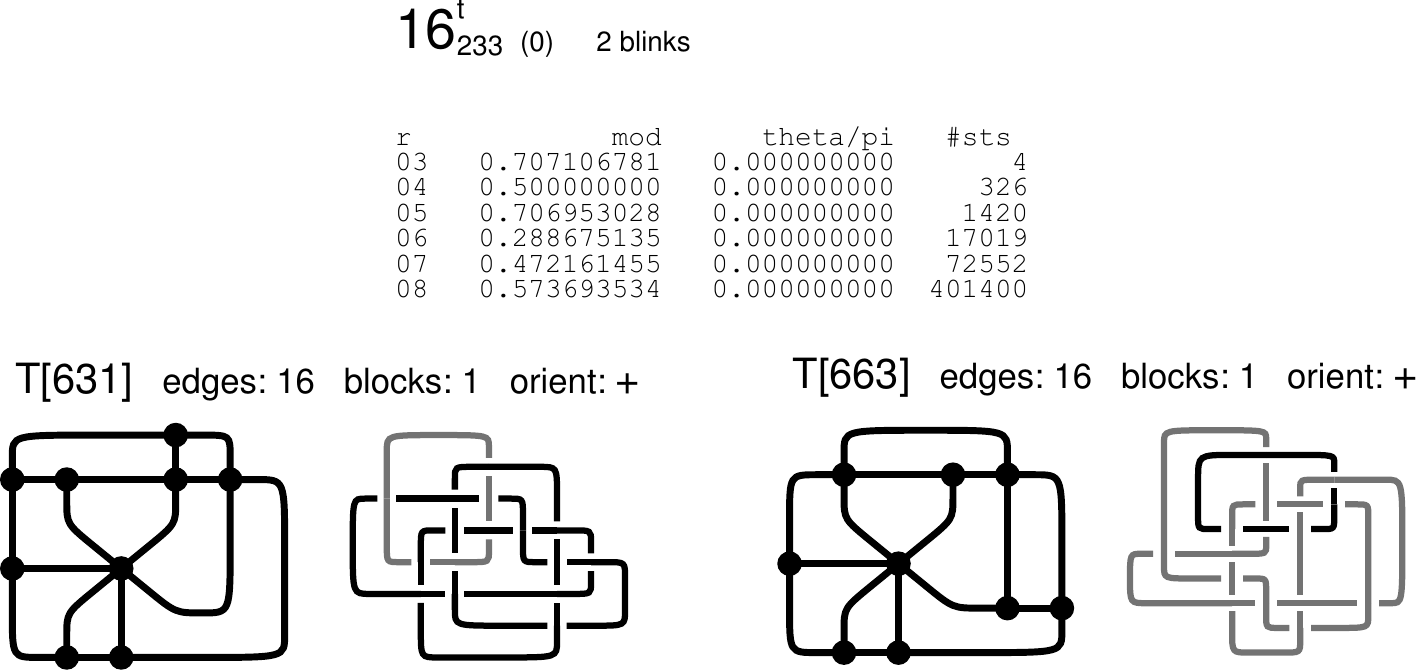}
\caption{\sf I do not know whether the above two manifolds are homeomorphic or not.}
\label{fig:T16_233}
\end{center}
\end{figure}

\section{Concluding remarks}
The elegant drawings of blinks and blackboard framed links produced by BLINK are 
possible due the groundbreaking algorithm of R. Tamasia \cite{tamassia1987egg}.
Lauro could implement the drawings very fast because we had at hand the implementation of
network flow algorithms he had done for a project to solve {\em practical timetable (!) problems}.
This is an example of the unicity in Mathematics, advocated by L. Lovasz in his famous essay \cite{lovasz1998om}.
To get the drawings one has to apply three times the full strength of network flow theory.
The drawings BLINK presents are in an integer grid and 
deterministically minimize the number of $\pi/2$-bents in the blackboarded framed links.
In particular, it permit us to deal with the unavoidable curls which adjust the integer framings in
the best possible way: we do not care about them. 
The drawings for the companion blinks require a slight modification: it replaces each $p$-valent vertex $p>4$,
by a $p$-polygon inducing 3-valent ones. The final result is massaged a bit to
produce aesthetically pleasing and unambiguous drawings.

As of this writing, C. Hodgson sent me some puzzling 
information (computed with a stronger version of SnapPy) about the first pair of 
manifolds. These are induced by $T[71]$ and $T[79]$, forming the $HG8QI$-class $14_{24}^t$.
They are non-homeomorphic 3-manifolds as first shown by N. Dunfield.
They are homology $\mathbb{Z}$-spheres which have the same WRT-invariants 
(according to BLINK), and quoting Craig {\em ``the same volume (around 24.8) and the same lenght spectra (up to 12 decimals):
the (complex) length of the first geodesic of $T[71]$ is 
0.4749346632398791 + 0i  (of multiplicity 1) 
and that of $T[79]$ is 0.4749346632399361 + 0 i  (of multiplicity 1).''}

\noindent
Here are the Dirichlet domains:
\begin{center}
\includegraphics[width=7cm]{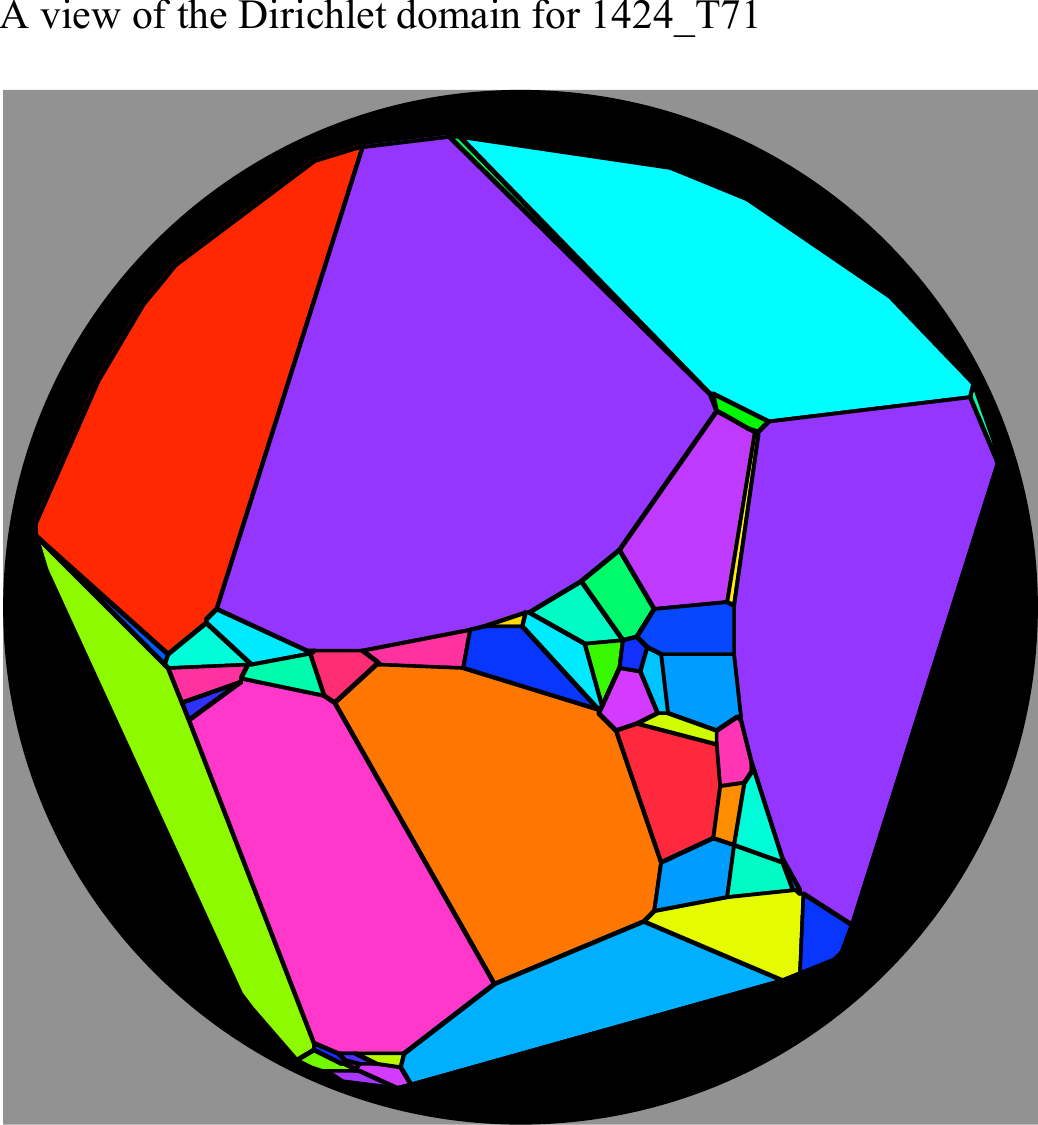}
\includegraphics[width=7cm]{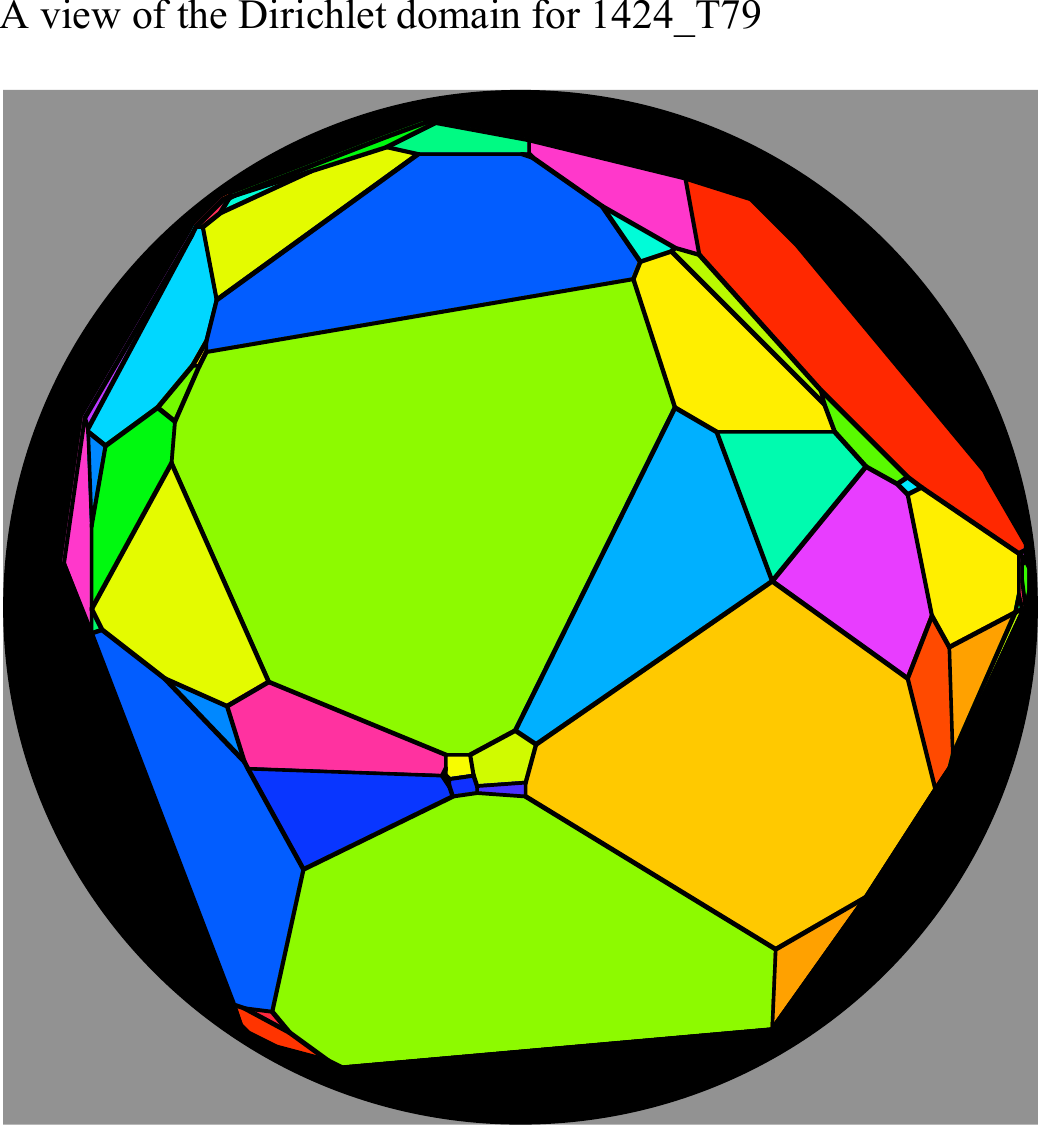}
\end{center}

\noindent
Craig found another proof that $T[71]$ and $T[79]$ are non-homemorphic:
{\em
``Now we can drill out the shortest geodesics using SnapPea to obtain
one-cusped manifolds (the manifold files are attached).
Then SnapPea's isometry checker (which uses the canonical cell decompositions) 
shows that these cusped manifolds are not isometric. Hence the original
closed manifolds are not homeomorphic.
This gives another proof that  $T[71]$ and $T[79]$ are distinct!''}

\noindent
{\bf A final challenge:} 
From the data I could get so far, if a closed orientable 3-manifold is hyperbolic, 
it seems that the WRT-invariants determine its volume. 
Prove it or disproved it.

\bibliographystyle{plain}
\bibliography{bibtexIndex.bib}

\begin{thebibliography}{10}

\bibitem{snappy}
M.~Culler, N.~M. Dunfield, and J.~Weeks.
\newblock Snappy, a computer program for studying the geometry and topology of
  3-manifolds, 2012.

\bibitem{gap2002gap}
GAP Group et~al.
\newblock Gap -- {G}roups, {A}lgorithms and {P}rogramming, version 4.3, 2002.

\bibitem{hodgson1994symmetries}
C.~D. Hodgson and J.~R. Weeks.
\newblock Symmetries, isometries and length spectra of closed hyperbolic
  three-manifolds.
\newblock {\em Experimental Mathematics}, 3(4):261--274, 1994.

\bibitem{kauffman1994tlr}
L.H. Kauffman and S.~Lins.
\newblock {Temperley-Lieb Recoupling Theory and Invariants of 3-manifolds}.
\newblock {\em Annals of Mathematical Studies, Princeton University Press},
  134:1--296, 1994.

\bibitem{Klarreich1012}
E~Klarreich.
\newblock Getting into shapes: from hyperbolic geometry to cube complexes.
\newblock {\em Simons Foundation}, October, 2012.

\bibitem{lins2007blink}
L.D. Lins.
\newblock Blink: a language to view, recognize, classify and manipulate
  3{D}-spaces.
\newblock {\em Arxiv preprint math/0702057}, 2007.

\bibitem{lins1995gca}
S.~Lins.
\newblock {\em {Gems, Computers, and Attractors for 3-Manifolds}}.
\newblock World Scientific, 1995.

\bibitem{linslinschallenge2013}
S.~Lins and S.~Lins.
\newblock A challenge to 3-manifold topologists and group algebraists.
\newblock {\em arXiv:1304.5964v4 [math.GT]}, 2013.

\bibitem{linsmachadoA2012}
S.~Lins and R.~Machado.
\newblock Framed link presentations of 3-manifolds by an ${O}(n^2)$ algorithm,
  {I}: gems and their duals.
\newblock {\em arXiv:1211.1953v2 [math.GT]}, 2012.

\bibitem{linsmachadoB2012}
S.~Lins and R.~Machado.
\newblock Framed link presentations of 3-manifolds by an ${O}(n^2)$ algorithm,
  {II}: colored complexes and boundings in their complexity.
\newblock {\em arXiv:1212.0826v2 [math.GT]}, 2012.

\bibitem{linsmachadoC2012}
S.~Lins and R.~Machado.
\newblock Framed link presentations of 3-manifolds by an ${O}(n^2)$ algorithm,
  {III}: geometric complex $\mathcal{H}_n^\star$ embedded into $\mathbb{R}^3$.
\newblock {\em arXiv:1212.0827v2 [math.GT]}, 2012.

\bibitem{lovasz1998om}
L.~Lovasz.
\newblock {One Mathematics}.
\newblock {\em The Berliner Intelligencer, Berlin}, pages 10--15, 1998.

\bibitem{rolfsen2003knots}
D.~Rolfsen.
\newblock {\em Knots and links}.
\newblock American Mathematical Society, 2003.

\bibitem{sage2012}
W.~Stein.
\newblock {\em {Sage}: {O}pen {S}ource {M}athematical {S}oftware ({V}ersion
  2.10.2)}.
\newblock The Sage~Group, 2008.
\newblock {\tt http://www.sagemath.org}.

\bibitem{stillwell1993classical}
J.~Stillwell.
\newblock {\em {Classical Topology and Combinatorial Group Theory}}.
\newblock Springer Verlag, 1993.

\bibitem{tamassia1987egg}
R.~Tamassia.
\newblock {On Embedding a Graph in the Grid with the Minimum Number of Bends}.
\newblock {\em SIAM Journal on Computing}, 16:421, 1987.

\bibitem{whitehead1997}
A.~N. Whitehead.
\newblock {\em The classical groups: their invariants and representations}.
\newblock Princeton {U}niversity {P}ress, 1997.

\end{thebibliography}

\vspace{5mm}
\begin{center}
\hspace{7mm}
\begin{tabular}{l}
   S\'ostenes L. Lins\\
   Centro de Inform\'atica, UFPE \\
   Av. Jornalista Anibal Fernandes s/n\\
   Recife, PE 50740-560 \\
   Brazil\\
   sostenes@cin.ufpe.br
\end{tabular}
\end{center}

\section{Appendix: text file to import
all links of the paper into SnapPy.
Prepared by Cristiana Nascimento}

\begin{verbatim}
% Link Projection 14-24_T71
2
   3    3
  15   15
16
   64    90
   59   232
  247   231
  247    91
  325   115
  208   116
  207   278
  366   275
  365   199
  162   204
  162   333
  274   333
  276   150
  108   149
  104   255
  327   253
16
   0    1
   1    2
   2    3
   3    0
   4    5
   5    6
   6    7
   7    8
   8    9
   9   10
  10   11
  11   12
  12   13
  13   14
  14   15
  15    4
14
   2    4
   5    1
   2    8
   8    5
   1    9
   6   11
   8   11
  12    2
   5   12
  13    1
  14    5
   9   14
  11   14
  15    8
-1

% Link Projection 14-24_T71
2
   3    3
  15   15
16
  235   205
  236   298
  363   296
  361   206
  192   120
  193   416
  279   416
  278    63
   56    62
   57   235
  311   234
  311   169
  111   173
  111   278
  406   272
  403   113
16
   0    1
   1    2
   2    3
   3    0
   4    5
   5    6
   6    7
   7    8
   8    9
   9   10
  10   11
  11   12
  12   13
  13   14
  14   15
  15    4
14
   1    6
   6    3
   0    9
   9    4
   9    6
   3   10
   4   11
  11    6
  12    9
  13    0
  13    2
   4   13
   6   13
   6   15
-1

% Link Projection  15-16_T118
2
   3    3
  17   17
18
  301   217
  304   350
  435   350
  428   213
  223   280
  457   271
  456   145
  335   146
  345   322
  137   323
  140   185
  410   175
  408    72
  177    79
  181   237
  379   240
  376   100
  225   110
18
   0    1
   1    2
   2    3
   3    0
   4    5
   5    6
   6    7
   7    8
   8    9
   9   10
  10   11
  11   12
  12   13
  13   14
  14   15
  15   16
  16   17
  17    4
15
   4    0
   4    2
   7    3
   7    4
   0    8
  10    7
  11    6
  10   13
   0   14
  14    7
   3   15
   6   15
  15   10
  17   10
  14   17
-1

% Link Projection 15-16_T119
2
   3    3
  17   17
18
  203    95
  205   183
  378   182
  378    93
  342   114
  448   110
  450   386
  248   386
  244   156
  422   158
  422   224
   97   222
   95   331
  285   328
  283    69
  156    72
  152   289
  345   288
18
   0    1
   1    2
   2    3
   3    0
   4    5
   5    6
   6    7
   7    8
   8    9
   9   10
  10   11
  11   12
  12   13
  13   14
  14   15
  15   16
  16   17
  17    4
15
   2    4
   1    7
   8    2
   7   10
   7   12
  13    1
  13    3
   8   13
  10   13
  10   15
  16    7
  13   16
   1   17
  17    8
  17   10
-1

% Link Projection 15-16_T181
2
   3    3
  17   17
18
   80    41
   84   250
  234   252
  228    50
  120   342
  118   213
  346   214
  341    81
  147    87
  152   392
  406   385
  404   123
  263   131
  267   301
   32   298
   28   170
  373   166
  378   338
18
   0    1
   1    2
   2    3
   3    0
   4    5
   5    6
   6    7
   7    8
   8    9
   9   10
  10   11
  11   12
  12   13
  13   14
  14   15
  15   16
  16   17
  17    4
15
   1    4
   2    5
   2    7
   8    1
   5    8
   6   11
   5   12
   4   13
  13    8
  15    0
  15    2
  15    6
   8   15
  12   15
   8   17
-1

% Link Projection 15-16_T205
2
   3    3
  17   17
18
  148   284
  146   384
  396   384
  391   284
  388   237
   95   241
   95   448
  262   446
  259   177
  453   175
  455   354
  183   357
  181   106
  322   102
  320   326
  215   325
  216   152
  392   146
18
   0    1
   1    2
   2    3
   3    0
   4    5
   5    6
   6    7
   7    8
   8    9
   9   10
  10   11
  11   12
  12   13
  13   14
  14   15
  15   16
  16   17
  17    4
15
   1    7
   7    3
   4    7
  10    2
   7   10
  11    3
   4   11
   3   13
  13    4
   8   13
  14    7
   3   15
  15    4
  13   16
  17    8
-1

% Link Projection 15-19_T122
2
   3    3
  17   17
18
  279   281
  276   378
  445   376
  445   280
  245   144
  247   422
  394   420
  388   234
  209   238
  207   175
  492   171
  490   338
  157   348
  162   210
  325   204
  331   477
  103   474
  100   155
18
   0    1
   1    2
   2    3
   3    0
   4    5
   5    6
   6    7
   7    8
   8    9
   9   10
  10   11
  11   12
  12   13
  13   14
  14   15
  15   16
  16   17
  17    4
15
   1    6
   3    6
   4    7
   4    9
   0   11
  11    2
  11    4
   6   11
  13    4
   8   13
  14    1
  14    3
   5   14
   7   14
  11   14
-1

% Link Projection 15-19_T148
2
   3    3
  17   17
18
  315   285
  318   422
  460   421
  457   285
  278   139
  279   450
  368   449
  363   204
  189   206
  193   354
  416   349
  415   244
  234   247
  235   174
  481   168
  483   380
  142   391
  143   143
18
   0    1
   1    2
   2    3
   3    0
   4    5
   5    6
   6    7
   7    8
   8    9
   9   10
  10   11
  11   12
  12   13
  13   14
  14   15
  15   16
  16   17
  17    4
15
   1    6
   6    3
   7    4
   0    9
   9    4
   9    6
   3   10
   4   11
  11    6
  12    7
   4   13
  15    0
  15    2
   4   15
   6   15
-1

% Link Projection 15-19_T188
2
   3    3
  17   17
18
   62    41
   62   155
  231   154
  231    43
   39    84
   42   232
  406   228
  401    54
  194    56
  198   304
  120   305
  114   116
  279   115
  281   273
  157   277
  156   195
  342   190
  343    87
18
   0    1
   1    2
   2    3
   3    0
   4    5
   5    6
   6    7
   7    8
   8    9
   9   10
  10   11
  11   12
  12   13
  13   14
  14   15
  15   16
  16   17
  17    4
15
   2    7
   8    1
   8    5
   1   10
  10    5
   2   11
  11    8
   5   12
  13    8
   5   14
  15    8
  12   15
  17    0
  17    2
   8   17
-1

% Link Projection 15-19_T208
2
   3    3
  17   17
18
   52    68
   53   205
  297   202
  296    70
   78    89
   80   263
  421   263
  419   120
  115   123
  115   299
  233   296
  232   155
  384   157
  385   226
  165   223
  167   343
  267   341
  267    91
18
   1    0
   2    1
   3    2
   0    3
   4    5
   5    6
   6    7
   7    8
   8    9
   9   10
  10   11
  11   12
  12   13
  13   14
  14   15
  15   16
  16   17
  17    4
15
   4    1
   2    7
   1    8
   8    5
  10    1
  10    5
  11    2
  13   10
   5   14
  14    9
   1   16
   5   16
   7   16
  16   11
  16   13
-1

% Link Projection 15-22_T128
2
   3    3
  17   17
18
  289   111
  289   274
  406   273
  406   109
  231   186
  441   188
  437    78
  169    78
  174   318
  472   320
  471   148
  344   149
  356   412
  108   409
  108   246
  430   237
  436   365
  233   367
18
   0    1
   1    2
   2    3
   3    0
   5    4
   6    5
   7    6
   8    7
   9    8
  10    9
  11   10
  12   11
  13   12
  14   13
  15   14
  16   15
  17   16
   4   17
15
   0    4
   2    4
  10    2
   5   10
   1   11
   4   11
  11    8
  14    0
  14    2
  14    7
  11   14
   8   15
  16   11
   8   17
  17   14
-1

% Link Projection 1522-T141
2
   3    3
  17   17
18
  110   219
  110   297
  280   299
  280   217
  171    89
   49    90
   51   348
  450   350
  445   121
  135   122
  135   261
  397   260
  398   387
  235   386
  236   185
  330   185
  333   417
  175   417
18
   0    1
   1    2
   2    3
   3    0
   5    4
   6    5
   7    6
   8    7
   9    8
  10    9
  11   10
  12   11
  13   12
  14   13
  15   14
  16   15
  17   16
   4   17
15
   9    3
   2   10
   6   11
  13    1
  13    3
   6   13
  10   13
  15    6
  10   15
  12   15
   1   17
   3   17
  17    6
  17    8
  17   10
-1

% Link Projection 15-22_T186
2
   3    3
  17   17
18
  236   239
  238   382
  325   382
  322   240
  272    60
  276   355
  418   355
  417   164
  201   163
  201    90
  459    87
  463   319
  152   318
  150   133
  384   134
  384   282
  186   280
  181    59
18
   0    1
   1    2
   2    3
   3    0
   4    5
   5    6
   6    7
   7    8
   8    9
   9   10
  10   11
  11   12
  12   13
  13   14
  14   15
  15   16
  16   17
  17    4
15
   3    4
   5    2
   4    7
   4    9
   0   11
   2   11
  11    4
  11    6
  13    4
   8   13
   7   14
  15    0
  15    2
   4   15
  13   16
-1

% Link Projection 15-22_T206
2
   3    3
  17   17
18
   86   265
   87   406
  407   406
  401   262
  119   143
  127   366
  320   366
  312    48
  399    49
  404   184
  235   182
  239   292
  363   292
  358    91
  180    92
  186   334
  445   335
  438   131
18
   0    1
   1    2
   2    3
   3    0
   4    5
   5    6
   6    7
   7    8
   8    9
   9   10
  10   11
  11   12
  12   13
  13   14
  14   15
  15   16
  16   17
  17    4
15
   3    4
   6    3
   9    6
   3   10
  11    6
   3   12
  12    9
  13    6
  14    3
  15    2
   6   15
   6   17
   8   17
  17   12
  17   14
-1

% Link Projection 16-42_T305
2
   3    3
  17   17
18
   94   114
   91   204
  295   208
  298   113
  335   135
   54   131
   55   226
  202   232
  204    94
  386    98
  389   285
  124   287
  121   178
  275   178
  273    63
  160    60
  156   264
  331   262
18
   0    1
   1    2
   2    3
   3    0
   4    5
   5    6
   6    7
   7    8
   8    9
   9   10
  10   11
  11   12
  12   13
  13   14
  14   15
  15   16
  16   17
  17    4
16
   4    0
   4    2
   1    7
   7    3
   4    7
   1   11
  11    6
   7   12
   3   13
  13    4
  13    8
  15    1
   3   15
  15    4
   6   15
  12   15
-1

% Link Projection 16-42_T337
2
   3    3
  17   17
18
  259   100
  259   258
  378   260
  377   105
  308   321
  304    60
  193    62
  196   223
  438   222
  444   392
  147   388
  142   140
  343   131
  346   320
  412   316
  406   171
  225   180
  230   316
18
   0    1
   1    2
   2    3
   3    0
   4    5
   5    6
   6    7
   7    8
   8    9
   9   10
  10   11
  11   12
  12   13
  13   14
  14   15
  15   16
  16   17
  17    4
16
   4    1
   4    3
   0    7
   7    2
   7    4
   0   11
  11    4
  11    6
   1   12
  12    7
  14    7
  15    0
   2   15
   4   15
  15   12
   7   16
-1

% Link Projection 16-42_T387
2
   3    3
  15   15
16
  224   129
  223   281
  366   295
  364   133
  303   352
  302    64
  169    62
  166   212
  438   215
  438    95
  261    96
  271   247
  415   255
  410   172
  128   165
  130   347
16
   0    1
   1    2
   2    3
   3    0
   4    5
   5    6
   6    7
   7    8
   8    9
   9   10
  10   11
  11   12
  12   13
  13   14
  14   15
  15    4
16
   1    4
   3    4
   7    0
   2    7
   7    4
   4    9
  10    3
  10    7
  11    2
   4   11
   7   12
   0   13
  13    2
   4   13
  13    6
  13   10
-1

% Link Projection 16-42_T420
2
   3    3
  17   17
18
  192   177
  195   280
  392   279
  388   175
  223   252
  226   338
  292   340
  281   131
  126   138
  126   217
  349   215
  341   152
  250   157
  260   313
  427   310
  422   110
  323   122
  328   256
18
   0    1
   1    2
   2    3
   3    0
   4    5
   5    6
   6    7
   7    8
   8    9
   9   10
  10   11
  11   12
  12   13
  13   14
  14   15
  15   16
  16   17
  17    4
16
   4    1
   6    1
   3    6
   0    9
   6    9
   3   10
   6   11
   1   12
  12    3
   9   12
  13    6
  16    3
   9   16
  11   16
  17    6
  12   17
-1

% Link Projection 16-56_T320
3
   3    3
  13   13
  17   17
18
  215   187
  216   346
  331   344
  324   185
  260   102
  264   316
  371   316
  368   270
  113   274
  111    56
  291    53
  291   295
  408   295
  402   101
  173   142
  372   143
  370   217
  167   215
18
   0    1
   1    2
   2    3
   3    0
   4    5
   5    6
   6    7
   7    8
   8    9
   9   10
  10   11
  11   12
  12   13
  13    4
  14   15
  15   16
  16   17
  17   14
16
   3    4
   5    2
   0    7
   7    2
   7    4
  10    3
  10    7
   2   11
  11    6
  10   13
   4   14
  14   10
  16    0
   2   16
   4   16
  16   10
-1

% Link Projection 16-56_T357
3
  11   11
  15   15
  19   19
20
   83    46
   84   226
  343   227
  345   303
  227   301
  229   171
  300   172
  302   335
  170   334
  170   129
  387   128
  389    42
  273    90
  267   269
  433   263
  429    95
  125    62
  130   194
  346   193
  342    66
20
   1    0
   2    1
   3    2
   4    3
   5    4
   6    5
   7    6
   8    7
   9    8
  10    9
  11   10
   0   11
  12   13
  13   14
  14   15
  15   12
  16   17
  17   18
  18   19
  19   16
16
   1    4
   1    6
   3    6
   8    1
  12    1
  12    5
   9   12
  13    2
   6   13
  10   15
   4   17
   6   17
  17    8
  17   12
  18    9
  15   18
-1

% Link Projection 16-56_T405
3
   3    3
   7    7
  17   17
18
   88    75
   89   200
  270   201
  269    78
  198   106
  199   297
  428   302
  428   103
  392   257
   38   265
   32   134
  313   134
  319   338
  135   347
  128   166
  238   162
  248   384
  388   387
18
   0    1
   1    2
   2    3
   3    0
   4    5
   5    6
   6    7
   7    4
   8    9
   9   10
  10   11
  11   12
  12   13
  13   14
  14   15
  15   16
  16   17
  17    8
16
   4    1
   2    7
   8    4
  10    0
  10    2
   4   10
  11    5
   8   11
   1   13
  13    8
  14    4
   1   15
   5   15
  15    8
  15   12
   5   17
-1

% Link Projection 16-56_T421
3
   3    3
  13   13
  17   17
18
   87   110
   91   187
  267   179
  261   105
  150   286
  140    69
   65    70
   69   243
  354   230
  364   355
  193   362
  183    64
  321    63
  333   275
  110   146
  116   325
  302   311
  285   134
18
   0    1
   1    2
   2    3
   3    0
   4    5
   5    6
   6    7
   7    8
   8    9
   9   10
  10   11
  11   12
  12   13
  13    4
  14   15
  15   16
  16   17
  17   14
16
   4    1
   4    3
   7    4
   1   10
   3   10
  10    7
  12    7
  13   10
   1   14
  14    7
  10   15
   7   16
  16   13
  17    2
  17    4
  10   17
-1

% Link Projection 16-140_T423
2
   3    3
  19   19
20
  132   101
  129   284
  225   286
  227    96
   78    67
   74   172
  339   170
  342   243
  162   251
  158   317
  307   311
  305    64
  431    64
  425   365
  192   365
  191   131
  388   130
  392   207
  262   204
  265    64
20
   0    1
   1    2
   2    3
   3    0
   4    5
   5    6
   6    7
   7    8
   8    9
   9   10
  10   11
  11   12
  12   13
  13   14
  14   15
  15   16
  16   17
  17   18
  18   19
  19    4
16
   0    5
   2    5
   7    2
   8    1
  10    5
  10    7
   1   14
   5   14
  14    7
  14    9
  15    2
  15   10
   6   17
  17   10
   5   18
  18   15
-1

% Link Projection 16-140_T444
2
   3    3
  17   17
18
  275   293
  274   372
  442   369
  441   287
  219   416
  342   414
  337   203
  141   206
  138   267
  407   263
  406   328
   85   337
   85   243
  296   237
  297   452
   43   450
   45   165
  225   160
18
   0    1
   1    2
   2    3
   3    0
   4    5
   5    6
   6    7
   7    8
   8    9
   9   10
  10   11
  11   12
  12   13
  13   14
  14   15
  15   16
  16   17
  17    4
16
   1    5
   3    5
   5    8
   9    3
   0   10
   5   10
   7   12
  13    1
  13    3
   4   13
   8   13
  10   13
  17    6
  17    8
  10   17
  12   17
-1

% Link Projection 16-140_T445
2
   3    3
  19   19
20
  177    77
  177   158
  409   158
  403    79
  124    49
  125   298
  463   297
  459    39
  336    38
  335   353
  207   351
  204   270
  380   271
  379   128
  247   129
  255   418
  159   416
  156   231
  286   227
  283    44
20
   0    1
   1    2
   2    3
   3    0
   5    4
   6    5
   7    6
   8    7
   9    8
  10    9
  11   10
  12   11
  13   12
  14   13
  15   14
  16   15
  17   16
  18   17
  19   18
   4   19
16
   1    8
   3    8
   5    8
   5   10
   8   11
   1   12
   8   13
   1   14
   5   14
   9   14
  11   14
   5   16
  14   17
   1   18
   3   18
  13   18
-1

% Link Projection 16-140_473
2
   3    3
  17   17
18
  103    72
  107   221
  233   228
  233    65
   82   163
  339   162
  347   300
  276   302
  270   126
  189   134
  203   375
  387   371
  372   258
  151   262
  141    99
  303    95
  315   333
   82   326
18
   0    1
   1    2
   2    3
   3    0
   4    5
   5    6
   6    7
   7    8
   8    9
   9   10
  10   11
  11   12
  12   13
  13   14
  14   15
  15   16
  16   17
  17    4
16
   4    0
   2    4
   4    7
   8    2
   9    1
   4    9
   5   12
   7   12
  12    9
   1   13
  13    4
   2   14
  15    4
  15    6
  12   15
   9   16
-1

% Link Projection 16-141_T424
2
   3    3
  19   19
20
  141    83
  147   191
  337   190
  326    79
  214    59
  228   332
  416   328
  400   105
  289   106
  307   301
   93   306
   90   242
  378   239
  369   136
  185   137
  192   280
  262   274
  256   207
   90   211
   82    61
20
   0    1
   1    2
   2    3
   3    0
   4    5
   5    6
   6    7
   7    8
   8    9
   9   10
  10   11
  11   12
  12   13
  13   14
  14   15
  15   16
  16   17
  17   18
  18   19
  19    4
16
   4    1
   4    3
   2    7
   1    8
   4    9
   4   11
   8   11
  13    2
  13    4
   8   13
   1   14
  11   14
  15    4
  11   16
  17    4
  14   17
-1

% Link Projection 16-141_T468
2
   3    3
  17   17
18
   59   290
   61   386
  216   387
  211   286
   33   256
   33   320
  354   313
  349   159
  269   157
  274   352
  183   357
  178   204
  407   201
  417   430
   86   428
   82   176
  316   177
  315   247
18
   0    1
   1    2
   2    3
   3    0
   4    5
   5    6
   6    7
   7    8
   8    9
   9   10
  10   11
  11   12
  12   13
  13   14
  14   15
  15   16
  16   17
  17    4
16
   0    5
   5    2
   8    5
   2    9
   3   10
  10    5
   6   11
   8   11
  14    1
  14    3
   5   14
  15    8
  11   16
  17    8
  10   17
  17   14
-1

% Link Projection 16-142_T425
2
  13   13
  17   17
18
   87    60
   86   244
  387   237
  383   165
  122   177
  119   364
  303   363
  301    89
  154    90
  150   393
  331   387
  336   199
  192   208
  190    60
  276   116
  271   260
  452   259
  451   114
18
   0    1
   1    2
   2    3
   3    4
   4    5
   5    6
   6    7
   7    8
   8    9
   9   10
  10   11
  11   12
  12   13
  13    0
  14   15
  15   16
  16   17
  17   14
16
   4    1
   1    6
   3    6
   1    8
   8    3
   8    5
  10    1
   6   11
   3   12
  12    7
  14    1
  14    3
  11   14
   6   15
  15   10
   6   17
-1

% Link Projection 16-142_T435
2
   3    3
  17   17
18
   90    77
   90   387
  253   387
  250    77
  127   101
  123   307
  439   305
  436   207
  158   217
  158   353
  482   354
  482   264
  331   270
  328   149
  200   148
  205   428
  370   426
  366   100
18
   0    1
   1    2
   2    3
   3    0
   4    5
   5    6
   6    7
   7    8
   8    9
   9   10
  10   11
  11   12
  12   13
  13   14
  14   15
  15   16
  16   17
  17    4
16
   5    2
   2    7
   5    8
   2    9
   6   11
   7   12
  13    2
  14    1
  14    5
   7   14
   9   14
  16    5
  16    7
   9   16
  11   16
   2   17
-1

% Link Projection 16-149_T437
2
   3    3
  19   19
20
   85    88
   91   341
  200   338
  195    89
   33   375
   29   267
  221   270
  216   115
  274   117
  282   303
  148   303
  148   183
  251   183
  243    88
  341    90
  347   223
  120   224
  119   147
  307   144
  317   373
20
   0    1
   1    2
   2    3
   3    0
   5    4
   6    5
   7    6
   8    7
   9    8
  10    9
  11   10
  12   11
  13   12
  14   13
  15   14
  16   15
  17   16
  18   17
  19   18
   4   19
16
   5    0
   5    2
   2    9
  10    5
  11    2
   6   11
   7   12
   2   15
  15    6
   8   15
  15   10
   2   17
  17    6
  17    8
  12   17
  15   18
-1

% Link Projection 16-149_T449
2
   3    3
  17   17
18
   71   199
   76   287
  242   290
  239   192
  153   159
  158   369
  364   368
  366   255
  194   252
  189   118
   26   127
   29   324
  312   322
  307   216
  112   224
  116   409
  285   405
  271   160
18
   0    1
   1    2
   2    3
   3    0
   5    4
   6    5
   7    6
   8    7
   9    8
  10    9
  11   10
  12   11
  13   12
  14   13
  15   14
  16   15
  17   16
   4   17
16
   4    1
   4    3
   2    7
   3    8
  11    4
  12    7
  13    2
  13    4
   8   13
   1   14
  14   11
   5   16
   7   16
  16   11
  16   13
   8   17
-1

% Link Projection 16-149_T464
2
   3    3
  19   19
20
   83    78
   83   373
  240   372
  236    82
  117   137
  118   274
  365   271
  365    99
  298   100
  296   247
  161   252
  165   182
  337   181
  339    76
  415    74
  413   211
  191   214
  193   317
  440   315
  442   124
20
   0    1
   1    2
   2    3
   3    0
   5    4
   6    5
   7    6
   8    7
   9    8
  10    9
  11   10
  12   11
  13   12
  14   13
  15   14
  16   15
  17   16
  18   17
  19   18
   4   19
16
   2    5
   9    2
  11    2
   8   11
   7   12
   2   15
   6   15
  15    8
   5   16
  16    9
  17    2
   2   19
  19    6
  19    8
  12   19
  14   19
-1

% Link Projection 16-149_T537
2
   3    3
  17   17
18
  196   113
  194   299
  306   301
  304   108
   53   187
  376   186
  377    83
  158    76
  161   389
  367   392
  369   251
  119   259
  116   443
  261   441
  259   146
  413   141
  415   349
   54   341
18
   0    1
   1    2
   2    3
   3    0
   4    5
   5    6
   6    7
   7    8
   8    9
   9   10
  10   11
  11   12
  12   13
  13   14
  14   15
  15   16
  16   17
  17    4
16
   0    4
   2    4
   4    7
  10    0
  10    2
   7   10
   1   13
   8   13
  13   10
   4   13
  14    2
   5   14
  16    7
  16    9
  11   16
  13   16
-1

% Link Projection 16-233_T631
2
   3    3
  17   17
18
   77    59
   74   193
  193   196
  189    63
   39   119
  235   115
  238   200
  376   201
  381   276
  154   272
  151    93
  289    92
  296   324
  120   322
  120   157
  344   160
  346   237
   39   241
18
   0    1
   1    2
   2    3
   3    0
   4    5
   5    6
   6    7
   7    8
   8    9
   9   10
  10   11
  11   12
  12   13
  13   14
  14   15
  15   16
  16   17
  17    4
16
   4    0
   4    2
   9    1
   9    4
   2   10
   6   11
   8   11
   1   13
   2   14
  14    5
  14    9
  11   14
  15    6
  16    9
  11   16
  13   16
-1

% Link Projection 16-233_T663
2
  15   15
  19   19
20
   51    94
  198    94
  201   255
   17   254
   18   318
  440   312
  438   147
  284   147
  286   419
  171   419
  170   181
  366   179
  365   366
  230   367
  231   283
   53   287
  104   121
  104   228
  324   227
  320   120
20
   0    1
   1    2
   2    3
   3    4
   4    5
   5    6
   6    7
   7    8
   8    9
   9   10
  10   11
  11   12
  12   13
  13   14
  14   15
  15    0
  16   17
  17   18
  18   19
  19   16
16
   7    4
   9    2
   9    4
  10    1
   7   10
   4   11
  12    7
   4   13
  14    9
   2   15
   1   17
  17    7
  17    9
  18    6
  10   18
   1   19
-1
\end{verbatim}

\end{document}